\documentclass[12pt,leqno]{article}

\usepackage[british]{babel}
\usepackage[tbtags]{amsmath}
\usepackage{amssymb,amsthm,amsopn,paralist,booktabs,url}


\theoremstyle{plain}
\newtheorem{theorem}{Theorem}[section]

\newtheorem{lemma}[theorem]{Lemma}

\theoremstyle{definition}

\theoremstyle{remark}

\newtheorem{remark}[theorem]{Remark}
\newtheorem*{acknowledgements}{Acknowledgements}

\numberwithin{equation}{section}
\numberwithin{table}{section}

\setdefaultenum{\upshape (i)}{}{}{}

\newenvironment{smallpars}{\footnotesize \setlength{\parindent}{0pt}
\setlength{\parskip}{10pt}}{\par}

\newcommand{\Lie}[1]{\textsl{#1}}
\newcommand{\lie}[1]{\mathfrak{#1}}
\DeclareMathOperator{\GL}{\Lie{GL}}

\DeclareMathOperator{\SO}{\Lie{SO}}
\DeclareMathOperator{\so}{\lie{so}}

\DeclareMathOperator{\SU}{\Lie{SU}}
\DeclareMathOperator{\su}{\lie{su}}
\DeclareMathOperator{\Un}{\Lie{U}}
\DeclareMathOperator{\un}{\lie{u}}

\newcommand{\inp}[2]{\langle #1, #2 \rangle}
\newcommand{\inpc}[2]{\inp{#1}{#2}_{\mathbb C}}

\DeclareMathOperator{\Ric}{Ric}

\newcommand{\lcf}{\lbrack\!\lbrack}
\newcommand{\rcf}{\rbrack\!\rbrack}
\newcommand{\real}[1]{\lcf #1 \rcf}
\newcommand{\alt}{\mathbf a}
\newcommand{\talt}{\tilde{\mathbf a}}

\newcommand{\Nt}{\widetilde\nabla}
\newcommand{\Nb}{\overline\nabla}

\newcommand{\Ieta}{\hat\eta}

\newcommand{\FCur}{\mathcal R}
\newcommand{\Kah}{\mathcal K}
\newcommand{\SRic}{\textsl{Ric}}
\DeclareMathOperator{\Cur}{\mathcal K}

\newcommand{\Kc}[1]{\mathcal K_{#1}}
\newcommand{\Cc}[1]{\mathcal C_{#1}}
\newcommand{\Wc}[1]{\mathcal W_{#1}}

\newcommand{\hook}{\lrcorner}
\newcommand{\Hodge}{\mathord{*}}
\newcommand{\T}{\checkmark}

\begin{document}
\title{Curvature of (Special) Almost Hermitian Manifolds}

\author{Francisco Mart\'\i n Cabrera and Andrew Swann}

\date{}

\maketitle

\begin{smallpars}
  \textbf{Abstract.} We study the curvature of almost Hermitian manifolds
  and their special analogues via intrinsic torsion and representation
  theory.  By deriving different foruml\ae\ for the skew-symmetric part of
  the $*$-Ricci curvature, we find that some of these contributions are
  dependent on the approach used, and for the almost Hermitian case we
  obtain tables that differ from those of Falcitelli, Farinola \& Salamon.
  We show how the exterior algebra may used to explain some of these
  variations.
  
  
\end{smallpars}

\section{Introduction}

In~\cite{Tricerri-Vanhecke:aH}, Tricerri and Vanhecke gave a complete
decomposition of the Riemannian curvature tensor~$R$ of an almost Hermitian
manifold $(M,I,\inp\cdot\cdot )$ into irreducible $\Un(n)$-components.
These divide naturally into two groups, one forming the
space~$\Kah=\Cur(\un(n))$ of algebraic curvature tensors for a K\"ahler
manifold, and the other, $\Kah^\perp$, being its orthogonal complement.

In~\cite{Falcitelli-FS:aH}, Falcitelli et al.\ showed that the components
of~$R$ in~$\Kah^\perp$ are linearly determined by the covariant derivative
$\nabla\xi$, where $\nabla$ is the Levi-Civita connection and $\xi$ is the
intrinsic torsion of the $\Un(n)$-structure on~$M$.  Gray and Hervella
\cite{Gray-H:16} showed that in general dimensions $\xi$ may be split into
four components $\xi_1,\dots,\xi_4$ under the action of $\Un(n)$.  By using
the minimal $\Un(n)$-connection $\Nt=\nabla + \xi$ of~$M$, Falcitelli et
al.\ display some tables which show whether or not the tensors $\Nt\xi_i$
and $\xi_i\odot\xi_j$ contribute to the components of~$R$ in~$\Kah^\perp$.
This provides a unified approach to many of the curvature results obtained
by Gray~\cite{Gray:curvature}.

The present paper is motivated by the interest in extending the above
results to special almost Hermitian manifolds.  These are defined as almost
Hermitian manifolds $(M,I,\inp\cdot\cdot)$ equipped with a complex volume
form $\Psi = \psi_+ + i \psi_-$.  Equivalently they are manifolds with
structure group~$\SU(n)$.  A detailed study of the intrinsic torsion $\eta
+ \xi$ of such manifolds was made in~\cite{Cabrera:special}, extending
results of Chiossi and Salamon \cite{Chiossi-S:SU3-G2}.  Here $\xi$~is the
intrinsic $\Un(n)$-torsion, as above, and $\eta$ is essentially a one-form.
There is much current interest in~$\SU(n)$-structures, partly as
generalisations of Calabi-Yau manifolds
\cite{Grantcharov-GP:CY-toric,Banos:MA6} and partly because of the r\^ole
played by torsion connections with holonomy~$\SU(n)$ in string
theory~\cite{Papadopoulos:brane,Gutowski-IP:calibrations}.

For $\SU(n)$ structures, the algebraic curvature tensors lie
in~$\Cur(\su(n))$ and are automatically Ricci-flat.  Therefore, one may
compute the Ricci curvature $\Ric$, and indeed the $*$-Ricci
curvature~$\Ric^*$, in terms of the intrinsic $\SU(n)$-torsion~$\eta+\xi$.
This enables us to find information about those $\SU(n)$-components of the
Riemannian curvature~$R$ which are determined by the tensors $\Ric$
and~$\Ric^*$.  Some of these components are contained in~$\Kah^\perp$ and
others are contained in~$\Kah$.  This will allow us, on the one hand, to
get more concrete information about some components of~$R$ contained
in~$\Kah^\perp$ and, on the other hand, to enlarge the tables of Falcitelli
et al.\ with columns related with some components contained in~$\Kah$.

In working out these contributions, we arrived at various alternative
formul\ae\ for certain curvature components purely in terms of the
intrinsic $\Un(n)$-torsion~$\xi$.  This leads to some entries in the tables
that are different from those obtained by Falcitelli et al.  To try to
account for this, we consider the identity $d^2=0$ in the exterior algebra.
Applying this to the K\"ahler $2$-form~$\omega$ and considering a
particular component indeed leads to a non-trivial relation between the
tensors contributing to the curvature.  One may view the relation
$d^2\omega=0$ as one way of taking account of some of the information that
the Levi-Civita connection $\nabla = \Nt - \xi$ is torsion-free.

The paper is organised as follows.  In~\S\ref{sec:preliminaries} we present
some preliminary material: definitions, results, notation, etc.  Then in
\S\ref{sec:curvature}, we derive some formul\ae\ relating curvature and
intrinsic torsion.  As an immediate application, we give an alternative
proof of Gray's result \cite{Gray:nearly-Kaehler} that any nearly K\"ahler
manifold of dimension six which is not K\"ahler is an Einstein manifold.
We then proceed to computing the contributions of different components of
the intrinsic torsion and its covariant derivative to the Ricci, $*$-Ricci
and Riemannian curvatures.  Because of the representation theory, this
behaves differently in dimensions $4$ and $6$ than in higher dimensions: in
dimension~$6$, $\xi$~splits into more $\SU(3)$-components; in
dimension~$4$, the space of curvature tensors is decomposed more finely
under the action of $\SU(2)$.  This motivates us to display results and
tables in two separate sections: \S\ref{sec:high}~for high dimensions, $2n
\geqslant 8$, and \S\ref{sec:low}~for dimensions six and four.  Finally, in
\S\ref{sec:exterior} we discuss identities derived from the exterior
algebra.

We remark that in this paper we will often use decompositions of tensor
products without providing explicit details, since such information can be
readily obtained via available computer programs.

\begin{acknowledgements}
  This work is supported by a grant from the MEC (Spain), project
  MTM2004-2644.  Andrew Swann thanks the Department of Fundamental
  Mathematics at the University of La Laguna for kind hospitality during
  the initial stages of this work.  Francisco Mart\'\i n Cabrera wishes to
  thank the Deparment of Mathematics and Computer Science at the University
  of Southern Denmark for kind hospitality whilst working on this project.
\end{acknowledgements}

\section{Preliminaries }
\label{sec:preliminaries}

An \emph{almost Hermitian} manifold is a $2n$-dimensional manifold~$M$, $n
>0$, with a $\Un(n)$-structure.  This means that $M$ is equipped with a
Riemannian metric $\inp\cdot\cdot$ and an orthogonal almost complex
structure~$I$.  Each fibre~$T_m M$ of the tangent bundle can be considered
as a complex vector space by defining $i x = Ix$.  We will write $T_m
M_{\mathbb C}$ when we are regarding~$T_m M$ as such a space.

We define a Hermitian scalar product $\inpc\cdot\cdot= \inp\cdot\cdot + i
\omega(\cdot,\cdot)$, where $\omega$ is the K\"ahler form given by $\omega
(x, y) = \inp x{Iy}$.  The real tangent bundle~$TM$ is identified with the
cotangent bundle~$T^*M$ by the map $x \mapsto \inp\cdot x=x$.  Analogously,
the conjugate complex vector space $\overline{T_m M_{\mathbb C}}$ is
identified with the dual complex space $T^*_m M_{\mathbb C}$ by the map $x
\mapsto \inpc\cdot x = x_{\mathbb C}$.  It follows immediately that
$x_{\mathbb C} = x + i Ix$.

If we consider the spaces $\Lambda^p T^*_m M_{\mathbb C}$ of skew-symmetric
complex forms, one can check that $x_{\mathbb C} \wedge y_{\mathbb C} = (x
+ i Ix) \wedge (y + i Iy)$.  There are natural extensions of the scalar
products $\inp\cdot\cdot$ and $\inpc\cdot\cdot$ to $\Lambda^p T^*_m M$ and
$\Lambda^p T^*_m M_{\mathbb C}$, defined respectively by
\begin{gather*}
  \inp ab = \frac1{p!} \sum_{i_1,\dots,i_p=1}^{2n} a(e_{i_1}, \dots,
  e_{i_p}) b(e_{i_1},\dots, e_{i_p}), \\
  \inpc{a_{\mathbb C}}{b_{\mathbb C}} = \frac1{p!} \sum_{i_1,\dots,i_p=1}^n
  a_{\mathbb C} (u_{i_1},\dots,u_{i_p}) \overline{b_{\mathbb
  C}(u_{i_1},\dots, u_{i_p})},
\end{gather*}
where $e_1, \dots , e_{2n}$ is an orthonormal basis for real vectors and
$u_1, \dots, u_n$ is a unitary basis for complex vectors.

The following conventions will be used in this paper.  If $b$ is a
$(0,s)$-tensor, we write
\begin{gather*}
  I_{(i)}b(X_1, \dots, X_i, \dots , X_s) = - b(X_1, \dots , IX_i, \dots ,
  X_s),\\
  I b(X_1,\dots,X_s) = (-1)^sb(IX_1,\dots,IX_s).
\end{gather*}

In~\cite{Tricerri-Vanhecke:aH}, Tricerri and Vanhecke gave a complete
decomposition of the Riemannian curvature tensor~$R$ of an almost Hermitian
manifold $(M,I,\inp\cdot\cdot)$ into irreducible $\Un(n)$-components.  As
was indicated above, some of these components, constituting a
$\Un(n)$-space denoted by~$\Kah=\Cur(\un(n))$, are the only components
which can occur when $M$ is a K\"ahler manifold.  In this text we will
follow the notation used in~\cite{Falcitelli-FS:aH} for such components.
Likewise, we will adopt the formalism used in \cite{Salamon:holonomy}
and~\cite{Falcitelli-FS:aH} for irreducible $\Un(n)$-modules.  Thus, for $n
\geqslant 2$,
\begin{equation*}
  \Kah = \Cc3 + \Kc1 + \Kc2,
\end{equation*}
where $\Cc3 \cong [\sigma_0^{2,2}]$, $\Kc1 \cong \mathbb R$, $\Kc2 \cong
[\lambda_0^{1,1}]$ and $+$~denotes direct sum.  We recall that
$\lambda_0^{p,q}$ is a complex irreducible $\Un(n)$-module coming from the
$(p,q)$-part of the complex exterior algebra and its corresponding dominant
weight in standard coordinates is given by $(1, \dots , 1,0, \dots , 0 , -1
, \dots , -1)$, where $1$ and $-1$ are repeated $p$ and $q$ times
respectively.

By analogy with the exterior algebra, there are also irreducible
$\Un(n)$-modules~$\sigma_0^{p,q}$ with dominant weights $(p,0, \dots , 0 ,
-q)$ coming from the symmetric algebra.  The notation~$\real V$ means the
real vector space underlying a complex vector space~$V$ and $[W]$~denotes a
real vector space which admits~$W$ as its complexification.

Moreover, let $\Ric$ and $\Ric^*$ respectively be the Ricci and $*$-Ricci
curvatures which are defined by
\begin{equation*}
  \Ric(X,Y) = \inp{R_{X,e_i} Y}{e_i}, \qquad
  \Ric^*(X,Y) = \inp{R_{X,e_i} IY}{Ie_i},
\end{equation*}
where $R_{X,Y} = \nabla_{[X,Y]} - [\nabla_X, \nabla_Y]$ and the summation
convention is used.

The components of the curvature~$R$ in $\Kc1$ and $\Kc2$ are determined by
the trace and the trace-free components of $\Ric_H + 3 \Ric_H^*$
respectively (see~\cite{Tricerri-Vanhecke:aH}), where $b_H$ indicates the
Hermitian part of a bilinear form~$b$, i.e., the part satisfying
$b_H(IX,IY)=b_H(X,Y)$.  Note that $\Ric^*_H$ coincides with the symmetric
part of $\Ric^*$.

The remaining components of~$R$, not included in~$\Kah$, are contained in a
$\Un(n)$-space denoted by~$\Kah^\perp$.  For $n\geqslant 4$, one
has~\cite{Falcitelli-FS:aH}
\begin{equation*}
  \Kah^\perp = \Kc{-1} + \Kc{-2} + \Cc4 + \Cc5 + \Cc6 + \Cc7 +
  \Cc8, 
\end{equation*}
where $\Kc{-1} \cong \mathbb R$, $\Kc{-2} \cong [\lambda_0^{1,1}]$, $\Cc4
\cong [\lambda_0^{2,2}]$, $\Cc5 \cong \real U$, $\Cc6 \cong
\real{\lambda^{2,0}} $, $\Cc7 \cong \real V$, and $\Cc8 \cong
\real{\sigma^{2,0}}$.  The irreducible $\Un(n)$-modules $U$ and $V$ have
dominant weights $(2,2, 0, \dots , 0)$ and $(2,1,0, \dots,0,-1)$
respectively.  For $n=3$, the decomposition of $\Kah^\perp$ is formed by
the same summands but omitting~$\Cc4$.  Finally, when $n=2$ we have to omit
$\Kc{-2}$, $\Cc4$ and~$\Cc7$.

We are dealing with $G$-structures where $G$ is a subgroup of the linear
group $\GL(m, \mathbb R)$.  If $M$ possesses a $G$-structure, then there
always exists a $G$-connection defined on~$M$.  Moreover, if $(M^m
,\inp\cdot\cdot)$ is an orientable $m$-dimensional Riemannian manifold and
$G$~a closed and connected subgroup of~$\SO(m)$, then there exists a unique
metric $G$-connection~$\Nt$ such that $\xi_x = \Nt_x - \nabla_x$ takes its
values in~$\lie g^\perp$, where $\lie g^\perp$ denotes the orthogonal
complement in~$\so(m)$ of the Lie algebra~$\lie g$ of~$G$ and $\nabla$~is
the Levi-Civita connection~\cite{Salamon:holonomy,Cleyton-S:intrinsic}.
The tensor~$\xi$ is the \emph{intrinsic torsion} of the $G$-structure and
$\Nt$ is called the \emph{minimal $G$-connection}.

For $\Un(n)$-structures, the minimal $\Un(n)$-connection is given by $\Nt =
\nabla + \xi$, with
\begin{equation} \label{torsion:xi}
  \xi_X Y = - \tfrac12 I\left( \nabla_X I \right) Y,
\end{equation}
see~\cite{Falcitelli-FS:aH}.  Since $\Un(n)$~stabilises the K\"ahler
form~$\omega$, it follows that $\Nt\omega = 0$.  Moreover, the equation
$\xi_X(IY) + I(\xi_XY) =0$ implies $\nabla\omega = - \xi\omega \in T^* M
\otimes \un(n)^\perp$.  Thus, one can identify the $\Un(n)$-components of
$\xi$ with the $\Un(n)$-components of $\nabla\omega$:
\begin{enumerate}
\item if $n=1$, $ \xi \in T^* M \otimes \un(1)^\perp = \{ 0 \}$;
\item if $n=2$, $ \xi \in T^* M \otimes \un(2)^\perp = \Wc2 +
  \Wc4$;
\item if $n \geqslant 3$, $ \xi \in T^* M \otimes \un(n)^\perp =
  \Wc1 + \Wc2 + \Wc3 + \Wc4$.
\end{enumerate}
Here the summands~$\Wc{i}$ are the irreducible $\Un(n)$-modules given by
Gray and Hervella in~\cite{Gray-H:16}, so $\Wc1 \cong
\real{\lambda^{3,0}}$, $\Wc2 \cong \real A$, $\Wc3 \cong
\real{\lambda^{2,1}_0}$ and $\Wc4 \cong \real{\lambda^{1,0}}$, where
$A\subset \lambda^{1,0}\otimes\lambda^{2,0}$~is the irreducible
$\Un(n)$-module with dominant weight $(2,1,0, \dots , 0 )$.  In the
following, $\xi_i$~will denote the component in~$\Wc{i}$ of the torsion
tensor~$\xi$.

In~\cite{Falcitelli-FS:aH}, Falcitelli et al.\ proved that the components
of~$R$ in~$\Kah^\perp$ are linearly determined by the covariant derivative
$\nabla\xi$ with respect to the Levi-Civita connection~$\nabla$.  To prove
this result, they consider the space $\FCur = \Kah + \Kah^\perp$ of
curvature tensors, we recall that $\FCur$ is the kernel of the mapping
$\odot^2 \left( \Lambda^2 T^*_m M \right) \to \Lambda^4 T^*_m M$ defined by
wedging two-forms together.  Then they deduce that the orthogonal
projection $\pi^\perp= \left( \pi_2 \circ \pi_1 \right)|_\FCur \colon \FCur
\to \Kah^\perp$ can be expressed as the restriction to~$\FCur$ of the
composition map $\pi_2 \circ \pi_1$, where $\pi_1 \colon \Lambda^2 T_m^*M
\otimes \Lambda^2 T_m^*M \to \Lambda^2 T_m^*M \otimes \un(n)^\perp$ is the
orthogonal projection and $\pi_2 \colon \Lambda^2 T_m^*M \otimes
\un(n)^\perp \to \Kah^\perp$ is a certain $\Un(n)$-equivariant
homomorphism.  Since we have the identity \cite{Falcitelli-FS:aH}
\begin{equation*}
  \begin{split}
    \pi_1(R)(X,Y,Z,W)
    &= \inp{\left( \nabla_X I \xi \right)_Y IZ}W
    - \inp{\left( \nabla_Y I \xi \right)_X IZ}W \\
    &= \inp{\left( \nabla_X \xi \right)_Y Z}W
    - \inp{\left( \nabla_Y \xi \right)_X Z}W \\
    &\qquad\qquad+ 2\inp{\xi_X\xi_YZ}W - 2\inp{\xi_Y\xi_XZ}W
  \end{split}
\end{equation*}
with the third and fourth summands in $\Lambda^2 T_m^*M \otimes 
\un(n)$, and $\pi_2$ is $\Un(n)$-equivariant, it follows that the
components of $\pi^\perp (R)$ in~$\Kah^\perp$ are linear functions
of the components of $\nabla\xi$.  Now, taking the
$\Un(n)$-connection $\Nt = \nabla + \xi$ into account, one obtains
\begin{equation} \label{ffsigualdad}
  \begin{split}
    \pi_1(R)(X,Y,Z,W)
    & = \inp{(\Nt_X\xi)_Y Z}W
    - \inp{(\Nt_Y\xi)_X Z}W \\
    &\quad + \inp{\xi_{\xi_X Y} Z}W - \inp{\xi_{\xi_Y X} Z}W.
  \end{split}
\end{equation}
From this equation and considering the image $\pi_2 \circ \pi_1 (R)$,
Falcitelli et al.\ give some tables which show whether or not the tensors
$\Nt\xi_i$ and $\xi_i \odot \xi_j$ contribute to the components of~$R$
in~$\Kah^\perp$.

Here we also consider manifolds equipped with an $\SU(n)$-structure.  Such
manifolds are called \emph{special almost Hermitian} manifolds.  They are
almost Hermitian manifolds $(M,I,\inp\cdot\cdot)$ equipped with a complex
volume form $\Psi = \psi_+ + i \psi_-$ such that $\inpc\Psi\Psi = 1$.  Note
that $I_{(i)} \psi_+ = \psi_-$.  See~\cite{Cabrera:special} for details and
more exhaustive information,
or~\cite{Bryant:splag,Joyce:holonomy,Hitchin:special-L}.

For a special almost Hermitian $2n$-manifold~$M$, we have the intrinsic
torsion $\eta + \xi \in T^* M \otimes \mathbb R\omega + T^*M \otimes
\un(n)^\perp = T^*M \otimes \su(n)^\perp$ and the minimal
$\SU(n)$-connection $\Nb = \nabla + \eta + \xi$.  Since $\Nb$ is metric and
$\eta \in T^* M \otimes \mathbb R\omega$, we have $\inp Y{\eta_X Z} =
\Ieta(X)\omega(Y,Z)$, where $\Ieta$~is a one-form.  Hence
\begin{equation*}
  \eta_X Y = \Ieta(X) IY.
\end{equation*}
In \cite{Cabrera:special} it is shown that the one-form~$\Ieta$ is given
by
\begin{equation*}
  -I\Ieta = \tfrac 1{2^{n-1}n} \Hodge(\Hodge d\psi_+\wedge\psi_+ + \Hodge
  d\psi_-\wedge\psi_-) - \tfrac1{2n} Id^*\omega,
\end{equation*}
where $\Hodge$ is the Hodge star operator and $d^*$ the coderivative.  This
formula simplifies for $n\geqslant3$ since then $\Hodge d\psi_+\wedge\psi_+
= \Hodge d\psi_-\wedge\psi_-$, and one sees that $nI\Ieta
-\frac12Id^*\omega$ is essentially the coefficient of~$\Psi$ in the
$(n,1)$-part of~$d\Psi$.  The other part of the intrinsic torsion $\xi \in
T^* M \otimes \un(n)^\perp$ is still given by equation~\eqref{torsion:xi}.

The tensors $\omega$, $\psi_+$ and $\psi_-$ are stabilised by the
$\SU(n)$-action, and $\Nb\omega = 0$, $\Nb\psi_+ = 0$ and $\Nb\psi_- = 0$.
Moreover, one can check $\eta\omega = 0$ and obtain $\nabla\omega =
-\xi\omega \in T^* M \otimes \un(n)^\perp$.  In general, the above
mentioned $\Un(n)$-spaces $\Wc i$ are also irreducible as $\SU(n)$-spaces.
The only exceptions are $\Wc1$ and $\Wc2$ when $n=3$.  In fact, for that
case, we have the following decompositions into irreducible
$\SU(3)$-components,
\begin{equation*}
  \Wc i = \Wc i^+ + \Wc i^-, \quad i=1,2,
\end{equation*}
where the space $\Wc i^+$ ($\Wc i^-$) consists of those tensors $a \in \Wc
i \subseteq T^*M \otimes \Lambda^2 T^* M$ such that the bilinear form
$r(a)$, defined by $ 2 r(a)(x,y) = \inp{x \hook \psi_+}{y \hook a}$, is
symmetric (skew-symmetric), see~\cite{Cabrera:special,Chiossi-S:SU3-G2}.
The components of the tensor~$\xi$ in $\Wc i^+$ and~$\Wc i^-$, $i=1,2$,
will be denoted by $\xi^+_i$ and $\xi^-_i$ respectively.  Writing $\eta \in
\Wc5 \cong T^*M$, the intrinsic $\SU(n)$-torsion $\xi + \eta$ is contained
in $(T^* M \otimes \un(n)^\perp) + \Wc5$.  The space $\Wc5$ is always
$\SU(n)$-irreducible.

From the equations $\Nb \psi_+ = 0$ and $\Nb \psi_- = 0$, we have $\nabla
\psi_+ = - \xi \psi_+ - \eta \psi_+$ and $\nabla \psi_- = - \xi \psi_- -
\eta \psi_-$.  Moreover, for $n\geqslant2$, in~\cite{Cabrera:special} it is
shown that
\begin{equation} \label{xipsi}
  \begin{gathered}
    \xi_X \psi_+, \xi_X \psi_- \in \real{ \lambda^{n-2,0} } \wedge
    \omega,\\
    \eta_X \psi_+ = n\, \Ieta(X) \psi_-, \qquad \eta_X \psi_- = - n\,
    \Ieta(X) \psi_+.
  \end{gathered}
\end{equation}

When considering curvature, note that the module $\Cc3=\Cur(\su(n))$
in~$\Kah$ consists of the algebraic curvature tensors for a metric with
holonomy algebra~$\su(n)$.

\section{Some curvature formul\ae}
\label{sec:curvature}

For special almost Hermitian $2n$-manifolds, results and tables given in
\cite{Falcitelli-FS:aH} are still valid with respect to the tensors
$\Nt\xi_i$ and $\xi_i \odot \xi_j$.  Here $\Nt =\Nb - \eta$ is the minimal
$\Un(n)$-connection, with $\Nb$ denoting the minimal $\SU(n)$-connection.

For $\SU(n)$-structures, the additional information coming from~$\eta$ will
allow us to compute the components of~$R$ in $\Kc1$ and $\Kc2$ in terms of
the intrinsic torsion $\eta + \xi$.  To achieve this, we compute the
difference between the Ricci and the $*$-Ricci curvatures.  In the first
instance we only need the almost Hermitian structure.

\begin{lemma} \label{ricmenosric}
  Let $M$ be an almost Hermitian $2n$-manifold, $n\geqslant 2$, with
  minimal $\Un(n)$-connection $\Nt = \nabla + \xi$, then
  \begin{equation*}
    \begin{split}
      \Ric^* (X,Y) - \Ric (X,Y)
      & = 2 \inp{(\nabla_{e_i} I\xi)_X IY}{e_i} - 2 \inp{(\nabla_X
      I\xi)_{e_i} IY}{e_i}, \\
      & = 2 \inp{(\Nt_{e_i}\xi)_X Y}{e_i} - 2 \inp{(\Nt_X\xi)_{e_i} Y}{e_i} \\
      &\qquad + 2 \inp{\xi_{\xi_{e_i}X}Y}{e_i} - 2 \inp{\xi_{\xi_X e_i}Y}{e_i}.
    \end{split}
  \end{equation*}
\end{lemma}

\begin{proof}
  It is straightforward to check
  \begin{equation} \label{difric}
    \Ric^*(X,Y) - \Ric(X,Y) =  - ( R_{X, e_i}\omega)(IY,e_i).
  \end{equation}
  However the so-called Ricci formula \cite[p.~26]{Besse:Einstein} implies
  \begin{equation} \label{ricident} - (R_{X,e_i}\omega)(IY,e_i) =
    \talt(\nabla^2\omega)_{X,e_i}(IY,e_i),
  \end{equation}
  where $\talt \colon T^* M \otimes T^* M \otimes \Lambda^2 T^* M \to
  \Lambda^2 T^* M \otimes \Lambda^2 T^* M$ is the skewing mapping.
  
  The required identities follow from equations \eqref{difric} and
  \eqref{ricident}, taking into account $\Nt\omega = 0$.
\end{proof}

The components of~$R$ in $\Kc{-1}$ and $\Kc{-2}$ are determined by the
trace and the trace-free parts of $\Ric^*_H - \Ric_H $.  Similarly, the
$\Cc6$-component of~$R$ is determined by the skew-symmetric (or
anti-Hermitian) part $\Ric^*_{AH}$ of~$\Ric^*$.  Moreover, the
anti-Hermitian part $\Ric_{AH}$ of the Ricci curvature, which satisfies
$\Ric_{AH}(IX,IY)= -\Ric_{AH}(X,Y)$, determines the component of~$R$
in~$\Cc8$.  These assertions motivate the expressions contained in the next
lemma.

\begin{lemma}
  Let $M$ be an almost Hermitian $2n$-manifold, $n\geqslant 2$, with
  minimal $\Un(n)$-connection $\Nt = \nabla + \xi$, then
  \begin{gather}
    \label{richmenosrich}
    \allowdisplaybreaks
    \begin{split}
      (\Ric^*_H - \Ric_H)(X,Y)
      &= \inp{(\Nt_{e_i}\xi)_X Y}{e_i}
      - \inp{(\Nt_X\xi)_{e_i} Y}{e_i}\\
      &\qquad
      + \inp{(\Nt_{e_i}\xi)_{IX} IY}{e_i} 
      - \inp{(\Nt_{IX}\xi)_{e_i} IY}{e_i}\\
      &\qquad
      + \inp{\xi_{\xi_{e_i}X}Y}{e_i}
      - \inp{\xi_{\xi_X e_i}Y}{e_i}\\
      &\qquad
      + \inp{\xi_{\xi_{e_i}IX}IY}{e_i}
      - \inp{\xi_{\xi_{IX} e_i} IY}{e_i},
    \end{split}
    \\
    \label{ricastsh}
    \begin{split}
      2\Ric^*_{AH}(X,Y)
      &= \inp{(\Nt_{e_i}\xi)_X Y}{e_i}
      - \inp{(\Nt_{e_i}\xi)_Y X}{e_i}\\
      &\qquad
      - \inp{(\Nt_{e_i}\xi)_{IX} IY}{e_i}
      + \inp{(\Nt_{e_i}\xi)_{IY} IX}{e_i}\\
      &\qquad
      - \inp{(\Nt_X\xi)_{e_i}Y}{e_i}
      + \inp{(\Nt_Y\xi)_{e_i}X}{e_i}\\
      &\qquad
      + \inp{(\Nt_{IX}\xi)_{e_i}IY}{e_i}
      - \inp{(\Nt_{IY}\xi)_{e_i}IX}{e_i}\\
      &\qquad
      + \inp{\xi_{\xi_X e_i}Y}{e_i}
      -  \inp{\xi_{\xi_Y e_i}X}{e_i}\\
      &\qquad
      - \inp{\xi_{\xi_{IX} e_i}IY}{e_i}
      +  \inp{\xi_{\xi_{IY} e_i}IX}{e_i},
    \end{split}
    \\
    \label{ricsh}
    \begin{split}
      2 \Ric_{AH}&(X,Y) \\
      =-&\inp{(\Nt_{e_i}\xi)_X Y}{e_i}
      + \inp{(\Nt_X\xi)_{e_i} Y}{e_i}
      - \inp{(\Nt_{e_i}\xi)_Y X}{e_i}\\
      &+ \inp{(\Nt_Y\xi)_{e_i} X}{e_i}
      + \inp{(\Nt_{e_i}\xi)_{IX} IY}{e_i}
      - \inp{(\Nt_{IX}\xi)_{e_i} IY}{e_i}\\
      &+ \inp{(\Nt_{e_i}\xi)_{IY} IX}{e_i}
      - \inp{(\Nt_{IY}\xi)_{e_i} IX}{e_i}
      - \inp{\xi_{\xi_{e_i}X}Y}{e_i}\\
      &+ \inp{\xi_{\xi_X e_i}Y}{e_i}
      - \inp{\xi_{\xi_{e_i}Y}X}{e_i}
      + \inp{\xi_{\xi_Y e_i}X}{e_i}
      + \inp{\xi_{\xi_{e_i}IX}IY}{e_i}\\
      &- \inp{\xi_{\xi_{IX} e_i}IY}{e_i}
      + \inp{\xi_{\xi_{e_i}IY}IX}{e_i}
      - \inp{\xi_{\xi_{IY}e_i}IX}{e_i}.
    \end{split}
  \end{gather}
\end{lemma}

\begin{proof}
  This follows directly from Lemma~\ref{ricmenosric} together with
  $\inp{\xi_{\xi_{e_i}X}Y}{e_i} = \inp{\xi_{\xi_{e_i}Y}X}{e_i}$.
\end{proof}

Up to this point, we have not said anything special about
$\SU(n)$-\hspace{0pt}structures.  We now give a first result that uses the
complex volume form~$\Psi$.

\begin{lemma}
  Let $M$ be a special almost Hermitian $2n$-manifold, $n\geqslant 2$, with
  complex volume form $\Psi = \psi_+ + i\psi_-$ and minimal
  $\SU(n)$-connection $\Nb = \nabla + \eta + \xi = \Nt + \eta $, then
  \begin{gather}
    \label{ricasteta}
    \Ric^*(X,Y)  = - n\, d\Ieta(X,IY) - \inp{\xi_X e_i}{\xi_{IY} Ie_i}, \\
    \label{riceta}
    \begin{split}
      \Ric(X,Y)
      &= - n\, d\Ieta(X,IY)
      - \inp{\xi_X e_i}{\xi_{IY} Ie_i}
      - 2\inp{(\Nt_{e_i}\xi)_X Y}{e_i}\\
      &\qquad
      + 2\inp{(\Nt_X\xi)_{e_i}Y}{e_i}
      - 2 \inp{\xi_{\xi_{e_i}X}Y}{e_i}
      + 2\inp{\xi_{\xi_X e_i}Y}{e_i}.
  \end{split}
\end{gather}
\end{lemma}

\begin{proof}
  Start by noting that $\inp{R_{X,Y} \psi_+}{\psi_-} = - 2^{n-2}
  \inp{R_{X,Y} Ie_i}{e_i}$.  Now, by the first Bianchi identity, we have
  \begin{equation} \label{volumholric}
    \inp{R_{X,Y} \psi_+}{\psi_-}  =  - 2^{n-1}\Ric^*(X,IY).
  \end{equation}
  On the other hand, using the Ricci formula $- R_{X,Y} \psi_+ = \talt
  (\nabla^2 \psi_+)_{(X,Y)}$ and taking $\Nb = \nabla + \eta + \xi$ into
  account, we obtain
  \begin{equation*}
    \begin{split}
      - R_{X,Y}\psi_+
      &= n\, d\Ieta(X,Y)\psi_-
      + n\,\Ieta(X)(\xi_Y\psi_-)
      - n\,\Ieta(Y)(\xi_X\psi_-) \\
      &\qquad + Y \hook \left( \nabla_X (\xi \psi_+) \right)
      - X \hook \left( \nabla_Y (\xi \psi_+) \right).
    \end{split}
  \end{equation*}
  Using the inclusions of~\eqref{xipsi}, we have $\inp{\xi_X
  \psi_+}{\psi_-}=0$, $\inp{\xi_X \psi_-}{\psi_-} =0$ and $\inp{Y \hook
  \left( \nabla_X (\xi \psi_+) \right)}{\psi_-} = - \inp{\xi_X (\xi_Y
  \psi_+)}{\psi_-}$.  This gives the following identity
  \begin{equation} \label{volumdeta}
    \inp{R_{X,Y} \psi_+}{\psi_-} =  - n 2^{n-1} d\Ieta(X,Y)
    - 2^{n-1} \inp{\xi_X e_i}{\xi_Y I e_i}.
  \end{equation}
  Using equations \eqref{volumholric}, \eqref{volumdeta} and
  Lemma~\ref{ricmenosric}, we obtain the required identities for $\Ric^*$
  and $\Ric$.
\end{proof}

The following theorem is an immediate consequence of the above Lemma.

\begin{theorem}
  Let $M$ be a special almost Hermitian $2n$-manifold, $n\geqslant2$, that
  is K\"ahler.  Then $\Ric^*=\Ric$ and
  \begin{enumerate}
  \item if $d\Ieta = \lambda\omega$, for some $\lambda \in
    \mathbb R\setminus\{0\}$, then the manifold is Einstein, or
  \item if the one-form $\Ieta$ is closed, then the manifold
    is Ricci flat.\qed
  \end{enumerate}
\end{theorem}

In~\cite{Gray:nearly-Kaehler}, Gray proved that any nearly K\"ahler (type
$\Wc1$) connected six-manifold which is not K\"ahler is Einstein.  Here we
give an alternative proof.

\begin{theorem}[Gray~\cite{Gray:nearly-Kaehler}]
  \label{alterproof}
  Let $M$ be a special almost Hermitian connected six-manifold of type
  $\Wc1^+ + \Wc1^- + \Wc5$ which is not of type~$\Wc5$.  Then $M$ is an
  Einstein manifold such that $\Ric =5\Ric^* = 5 \alpha\,\inp\cdot\cdot$,
  where $\alpha = (w_1^+)^2 + (w_1^-)^2$ with $\nabla \omega = w_1^+ \psi_+
  + w_1^- \psi_-$.
\end{theorem}

\begin{proof}
  We already know that $\alpha = (w_1^+)^2 + (w_1^-)^2$ is a positive
  constant and the one-form $\Ieta$ is closed (see
  \cite[Theorem~3.7]{Cabrera:special}).  On the other hand, since $\nabla
  \omega = - \xi \omega$ and $\nabla \omega = w_1^+ \psi_+ + w_1^- \psi_-$,
  we have
  \begin{equation*}
    2 \inp{Y}{\xi_X Z} = w_1^- \psi_+(X,Y,Z) - w_1^+ \psi_-(X,Y,Z).
  \end{equation*}
  Therefore, using 
  \begin{gather*}
    \inp{X \hook \psi_+}{Y \hook \psi_+} = \inp{X \hook
    \psi_-}{Y \hook \psi_-} = 2 \inp XY,\\
    \inp{X \hook \psi_+}{Y \hook \psi_-} = -2\omega(X,Y),
  \end{gather*}
  we get
  \begin{equation}
    \inp{\xi_X e_i}{\xi_Y e_i} = \inp{e_j}{\xi_X e_i} \inp{e_j}{\xi_Y e_i}
    = \alpha\, \inp XY. 
  \end{equation}
  Moreover, since $\xi \in \Wc1^+ + \Wc1^-$ and $\Nt$ is a
  $\Un(3)$-connection, the $(0,3)$-tensors $\inp\cdot{\xi_{\cdot}\cdot}$
  and $\inp\cdot{(\Nt_X \xi)_{\cdot}\cdot}$ are skew symmetric
  \cite{Gray-H:16}.  Thus, from \eqref{riceta}, we get
  \begin{equation*}
    \Ric (X,Y) = 5\inp{\xi_X e_i}{\xi_Y e_i} = 5 \alpha\, \inp XY.
  \end{equation*}
  We recall that $\inp Y{\xi_{IX} IZ} = - \inp Y{\xi_{X} Z}$, for $\xi \in
  \Wc1$, and note that the contractions $\inp{(\Nt_X \xi)_{e_i} Y}{e_i}$ and
  $\inp{(\Nt_{e_i} \xi)_X Y}{e_i}$ both vanish.  In fact, the last term is
  a skew-symmetric two-form and the remaining summands in the expression
  for~$\Ric$ are symmetric.
\end{proof}

\begin{remark}
  Theorem~\ref{alterproof} can be extended to connected almost Hermitian
  six-manifolds which are nearly K\"ahler and but not K\"ahler.  In fact,
  one can define a complex volume form on an open neighbourhood~$U$ of a
  point where $\nabla \omega \neq 0$ by using the $(3,0)$-component of this
  tensor.  Then, $U$~is a special almost Hermitian six-manifold of type
  $\Wc1^+ + \Wc1^- + \Wc5$.  Therefore, $\Ric =5\Ric^* = 5 \alpha\,
  \inp\cdot\cdot$ on~$U$.  Since the manifold is connected, it follows
  $\Ric = 5 \alpha\, \inp\cdot\cdot$ everywhere.
\end{remark}

The expressions \eqref{ricasteta} and \eqref{riceta} for $\Ric^*$ and
$\Ric$ allow us to compute $3\Ric^*_H + \Ric_H$ and study the contributions
of the intrinsic torsion of the $\SU(n)$-structure to the components of $R$
in $\Kc1$ and~$\Kc2$.

\begin{lemma}
  Let $M$ be a special almost Hermitian $2n$-manifold, $n\geqslant 2$,
  with minimal $\SU(n)$-connection $\Nb = \nabla + \eta + \xi = \Nt +
  \eta$, then
  \begin{equation}
    \label{tresricmasric}
    \begin{split}
      (3\Ric^*_H &+ \Ric_H)(X,Y)\\
      &= - 2n\, d\Ieta(X,IY) + 2n\, d\Ieta(IX,Y) -
      \inp{(\Nt_{e_i}\xi)_XY}{e_i}\\
      &\qquad+\inp{(\Nt_X\xi)_{e_i}Y}{e_i} - \inp{(\Nt_{e_i}\xi)_{IX}IY}{e_i} +
      \inp{(\Nt_{IX}\xi)_{e_i}IY}{e_i} \\
      &\qquad - \inp{\xi_{\xi_{e_i}X}Y}{e_i} + \inp{\xi_{\xi_{X}e_i}Y}{e_i} -
      \inp{\xi_{\xi_{e_i}IX}IY}{e_i} \\
      &\qquad + \inp{\xi_{\xi_{IX}e_i}IY}{e_i} - 2\inp{\xi_Xe_i}{\xi_{IY}Ie_i} -
      2\inp{\xi_Ye_i}{\xi_{IX}Ie_i}.
      \qed
    \end{split} 
  \end{equation}
\end{lemma}

To end this section, let us note an alternative to
equation~\eqref{ricastsh} is given by the following:

\begin{lemma}
  Let $M$ be an almost Hermitian $2n$-manifold, $n\geqslant 2$, with
  minimal $\Un(n)$-connection $\Nt = \nabla + \xi$, then
  \begin{equation}
    \label{otraricsh}
    \begin{split}
      \Ric^*_{AH} (X,Y)
      & = \inp{(\Nt_{e_i}\xi)_{Ie_i}IX}Y - \inp{\xi_{I\xi_{e_i}e_i}IX}Y\\
      & = - \inp{(\nabla_{e_i}I\xi)_{Ie_i}X}Y.
    \end{split}
\end{equation}
\end{lemma}

\begin{proof}
  We have
  \begin{equation*}
    \begin{split}
      - 2 \Ric^* (X,IY) - 2 \Ric^* (IX,Y)
      & =  \inp{R_{e_i,Ie_i}X}Y - \inp{R_{e_i,Ie_i}IX}{IY} \\
      & =  4\inp{\Nt_{e_i}\xi_{Ie_i}X}Y - 4 \inp{\xi_{Ie_i}\Nt_{e_i}X}Y\\
      &\qquad - 4\inp{\xi_{\Nt_{e_i}Ie_i}X}Y + 4\inp{\xi_{\xi_{e_i}Ie_i}X}Y
    \end{split}
  \end{equation*}
  from which the Lemma follows.
\end{proof}

\section{High dimensions}
\label{sec:high}

In this section, we consider special almost Hermitian manifolds of
dimension higher than or equal to eight.  For such manifolds, the
decomposition into $\SU(n)$-irreducible modules of the space of curvature
tensors $\FCur$ is the same as that coming from the action of $\Un(n)$.
Thus,
\begin{equation*}
  \FCur = \Kah + \Kah^\perp = \Cc3 + \Kc1 + \Kc2 + \Kc{-1} +
  \Kc{-2} + \Cc4 + \Cc5 + \Cc6 + \Cc7 + \Cc8,
\end{equation*}
where all $\Kc{i}$ and $\Cc{j}$ are also $\SU(n)$-irreducible spaces.  Our
aim here is to see whether or not different components of the intrinsic
torsion of the $\SU(n)$-structure contribute to the components of the
curvature.

We start by studying such contributions to the $\SU(n)$-components of the
Ricci and $*$-Ricci curvatures.  For $n \geqslant 3$, the spaces $\SRic$
and $\SRic^*$ of such tensors admit the following decompositions into
$\SU(n)$-irreducible modules
\begin{equation*}
  \SRic = \mathbb R \inp\cdot\cdot + [ \lambda_0^{1,1} ] +
  \real{\sigma^{2,0}}, \qquad \SRic^* = \mathbb R \inp\cdot\cdot +
  [\lambda_0^{1,1} ] + \real{ \lambda^{2,0} }. 
\end{equation*}

Taking into account the symmetry properties and types of the
Gray-Hervella's components~$\xi_i$ of~$\xi$, we obtain the following
result.

\begin{theorem}
  \label{thm:Ric}
  Let $M$ be a special almost Hermitian $2n$-manifold, $2n\geqslant 8$,
  with minimal $\SU(n)$-connection $\Nb = \nabla + \eta + \xi = \Nt + \eta
  $.  The tensors $d\Ieta$, $\Nb\xi$ and $\xi_i \odot \xi_j$ contribute to
  the components of the $*$-Ricci curvature~$\Ric^*$ via
  equation~\eqref{ricasteta} and to the Ricci curvature~$\Ric$ via
  equation~\eqref{riceta} if and only if there is a tick in the
  corresponding place in Table~\ref{tab:Ric8}.
  \qed
\end{theorem}

\begin{table}[tp]
  \centering
  \begin{tabular}{lcccccc}
    \toprule
    &\multicolumn{3}{c}{$\Ric^*$
    \eqref{ricasteta}}&\multicolumn{3}{c}{$\Ric$ \eqref{riceta}}\\
    \cmidrule(lr){2-4}
    \cmidrule(lr){5-7}
    $2n \geqslant 8$&$\mathbb
    R$&$[\lambda^{1,1}_0]$&$\real{\lambda^{2,0}}$&$\mathbb
    R$&$[\lambda^{1,1}_0]$&$\real{\sigma^{2,0}}$\\
    \midrule
    $d\Ieta$                &\T&\T&\T&\T&\T&  \\
    \midrule                                                  
    $\Nb\xi_1$, $\eta\xi_1$ &  &  &  &  &  &  \\
    $\Nb\xi_2$, $\eta\xi_2$ &  &  &  &  &  &\T\\
    $\Nb\xi_3$, $\eta\xi_3$ &  &  &  &  &\T&  \\
    $\Nb\xi_4$, $\eta\xi_4$ &  &  &  &\T&\T&\T\\
    \midrule                                                  
    $\xi_1 \otimes \xi_1$   &\T&\T&  &\T&\T&  \\
    $\xi_2 \otimes \xi_2$   &\T&\T&  &\T&\T&  \\
    $\xi_3 \otimes \xi_3$   &\T&\T&  &\T&\T&\T\\
    $\xi_4 \otimes \xi_4$   &\T&\T&  &\T&\T&\T\\
    \midrule                                                  
    $\xi_1 \odot \xi_2$     &  &\T&  &  &\T&  \\
    $\xi_1 \odot \xi_3$     &  &  &\T&  &  &\T\\
    $\xi_1 \odot \xi_4$     &  &  &\T&  &  &  \\
    $\xi_2 \odot \xi_3$     &  &  &\T&  &  &\T\\
    $\xi_2 \odot \xi_4$     &  &  &\T&  &  &\T\\
    $\xi_3 \odot \xi_4$     &  &\T&  &  &\T&  \\
    \bottomrule
  \end{tabular}
  \caption{Ricci curvatures, $2n\geqslant8$}
  \label{tab:Ric8}
\end{table}

Using in addition $\inp{\xi_{\xi_{X} e_i}Y}{e_i} = - \inp{\xi_X
e_i}{\xi_{e_i}Y}$ we get part~(i) of the following theorem.  Part~(ii) is
proved in~\cite{Falcitelli-FS:aH}.

\begin{theorem}
  \label{thm:Curv}
  Let $M$ be a special almost Hermitian $2n$-manifold, $2n\geqslant 8$,
  with minimal $\SU(n)$-connection $\Nb = \nabla + \eta + \xi = \Nt + \eta
  $, then
  \begin{enumerate}
  \item Using equations \eqref{richmenosrich}, \eqref{ricastsh},
    \eqref{ricsh} and \eqref{tresricmasric}, each of the tensors
    $\Nb\xi_i$, $\eta\xi_i$ and $\xi_i \odot \xi_j$ contributes to the
    components of~$R$ in $\Kc{1}$, $\Kc{2}$, $\Kc{-1}$, $\Kc{-2}$, $\Cc6$
    and $\Cc8$ if and only if there is a tick in the corresponding place in
    Table~\ref{tab:Cur8}.
  \item Taking the image $\pi_2 \circ \pi_1 (R)$ into account, where $\pi_1
    (R)$ is given by equation~\eqref{ffsigualdad}, each of the tensors
    $\Nb\xi_i$, $\eta\xi_i$ and $\xi_i \odot \xi_j$ contributes to the
    components of $R$ in $\Cc4$, $\Cc5$ and $\Cc7$ if and only if there is
    a tick in the corresponding place in Table~\ref{tab:Cur8}.  \qed
  \end{enumerate}
\end{theorem}

\begin{table}[tp]
  \newcommand{\F}{$\T^*$}
  \centering
  \begin{tabular}{lcccccccccc}
    \toprule
    &\multicolumn{2}{c}{\eqref{tresricmasric}}
    &\multicolumn{2}{c}{\eqref{richmenosrich}}
    &\eqref{ricastsh}
    &\eqref{otraricsh}
    &\eqref{ricsh}
    &\multicolumn{3}{c}{\cite{Falcitelli-FS:aH}}\\
    \cmidrule(rl){2-3}
    \cmidrule(rl){4-5}
    \cmidrule(rl){6-6}
    \cmidrule(rl){7-7}
    \cmidrule(rl){8-8}
    \cmidrule(rl){9-11}
    $2n \geqslant 8$&$\Kc{1}$&$\Kc{2}$&$\Kc{-1}$&$\Kc{-2}$&$\Cc6$&$\Cc6
    $&$\Cc8$&$\Cc4$&$\Cc5$&$\Cc7$\\
    \midrule
    $d\Ieta$                &\T&\T&  &  &  &  &  &  &  &  \\
    \midrule
    $\Nb\xi_1$, $\eta\xi_1$ &  &  &  &  &\T&\T&  &  &  &  \\
    $\Nb\xi_2$, $\eta\xi_2$ &  &  &  &  &\T&\T&\T&  &\T&\T\\
    $\Nb\xi_3$, $\eta\xi_3$ &  &\T&  &\T&  &\T&  &\T&  &\T\\
    $\Nb\xi_4$, $\eta\xi_4$ &\T&\T&\T&\T&\T&\T&\T&  &  &  \\
    \midrule
    $\xi_1 \otimes \xi_1$   &\T&\T&\T&\T&  &  &  &\F&  &  \\
    $\xi_2 \otimes \xi_2$   &\T&\T&\T&\T&  &  &  &\T&  &  \\
    $\xi_3 \otimes \xi_3$   &\T&\T&  &  &  &  &\T&  &  &\T\\
    $\xi_4 \otimes \xi_4$   &\T&\T&  &  &  &  &\T&  &  &  \\
    \midrule
    $\xi_1 \odot \xi_2$     &  &\T&  &\T&  &  &  &\T&  &  \\
    $\xi_1 \odot \xi_3$     &  &  &  &  &\T&  &\T&  &\T&\T\\
    $\xi_1 \odot \xi_4$     &  &  &  &  &\T&\T&  &  &  &  \\
    $\xi_2 \odot \xi_3$     &  &  &  &  &\T&  &\T&  &\T&\T\\
    $\xi_2 \odot \xi_4$     &  &  &  &  &\T&\T&\T&  &\T&\T\\
    $\xi_3 \odot \xi_4$     &  &\T&  &  &  &\T&  &  &  &\T\\
    \bottomrule
    \multicolumn{11}{r}{\strut\footnotesize $^*$absent when $2n=8$}
  \end{tabular}
  \caption{Curvature complementary to $\Cc3=\Cur(\su(n))$, $2n\geqslant8$.} 
  \label{tab:Cur8}
\end{table}

For part~(i), we wish to emphasise that the columns for $\Kc{-1}$,
$\Kc{-2}$, $\Cc6$ and $\Cc8$ are obtained by a different method to that
in~\cite{Falcitelli-FS:aH} and that for $\Cc6$ this even leads to a
different result.  In particular, we claim that the tensors $\Nt\xi_3$ and
$\xi_3 \odot \xi_4$ do not contribute to the $\Cc6$-component of~$R$, but
that $\Nt\xi_1$ and $\eta\xi_1$ do.  Thus the contributions of the
different tensors to the distinct components of~$R$ depend on the choice of
the current expression that we use; different expressions may lead to
different behaviour in the contributions.  For the $\Cc6$-component of~$R$,
we get a third formula from equation~\eqref{otraricsh}, which we also list
in Table~\ref{tab:Cur8}.  A partial explanation for these different results
will be given in~\S\ref{sec:exterior}.  Note that the entries for $\Cc6$ in
Table~\ref{tab:Cur8} only involve the intrinsic $\Un(n)$-torsion.  The
$\real{\lambda^{2,0}}$-column of Table~\ref{tab:Ric8} provides yet another
description of the $\Cc6$-component using the $\SU(n)$-structure.

\section{Low dimensions}
\label{sec:low}

In this section we consider in turn special almost Hermitian manifolds of
dimension six and four.

\subsection{Six dimensions}

The decomposition of the space of curvature tensors $\FCur$ into
irreducible $\SU(3)$-modules has the same subspaces as that for~$\Un(3)$.
Thus,
\begin{equation*}
  \FCur = \Kah + \Kah^\perp = \Cc3 + \Kc1 + \Kc2 + \Kc{-1} +
  \Kc{-2} + \Cc5 + \Cc6 + \Cc7 + \Cc8,
\end{equation*}
with $\Kc{i}$ and $\Cc{j}$ all $\SU(3)$-irreducible.  As we noted above,
the summand $\Cc4$ is absent in this dimension.  On the other hand, the
$\Un(3)$-intrinsic torsion splits under $\SU(3)$ as $\xi = \xi_1^+ +
\xi_1^- + \xi_2^+ + \xi_2^- + \xi_3 + \xi_4$, where $\xi_i = \xi_i^+ +
\xi_i^-$, $i=1,2$.  This was briefly described in~\S\ref{sec:preliminaries}
and more detailed information is contained in \cite{Chiossi-S:SU3-G2}
and~\cite{Cabrera:special}.

The next result concerns the contributions of the components of $\xi$ to
the components of the Ricci and the $*$-Ricci curvatures and then to
the curvature components complementary to $\Cc3$.

\begin{theorem}
  Let $M$ be a special almost Hermitian $6$-manifold with
  $\SU(3)$-connection $\Nb = \nabla + \eta + \xi = \Nt + \eta $.  The
  tensors $d\Ieta$, $\Nb\zeta$, $\eta\zeta$ and $\zeta \odot \vartheta$,
  for $\zeta, \vartheta = \xi_1^+, \xi_1^-, \xi_2^+, \xi_2^-, \xi_3, \xi_4$
  contribute the components of $\Ric^*$ and $\Ric$ if and only if there is
  a tick in the corresponding place in Table~\ref{tab:Ric6}.
  
  The corresponding contributions to the curvature components $\Kc{1}$,
  $\Kc{2}$, $\Kah_{-1}$, $\Kc{-2}$, $\Cc6$ and $\Cc8$, via equations
  \eqref{richmenosrich}, \eqref{ricastsh}, \eqref{ricsh} and
  \eqref{tresricmasric}, and to the components $\Cc5$ and $\Cc7$ via
  $\pi_2\circ\pi_1(R)$ \cite{Falcitelli-FS:aH} are given in
  Table~\ref{tab:Cur6}. \qed
\end{theorem}

\begin{table}[tp]
  \centering
  \begin{tabular}{lcccccc}
    \toprule
    &\multicolumn{3}{c}{$\Ric^*$
    \eqref{ricasteta}}&\multicolumn{3}{c}{$\Ric$ \eqref{riceta}}\\
    \cmidrule(lr){2-4}
    \cmidrule(lr){5-7}
    $2n=6$&$\mathbb R$&$[\lambda^{1,1}_0]$&$\real{\lambda^{2,0}}$&$\mathbb
    R$&$[\lambda^{1,1}_0]$&$\real{\sigma^{2,0}}$\\ 
    \midrule
    $d\Ieta$                   &\T&\T&\T&\T&\T&  \\
    \midrule
    $\Nb\xi^\pm_1$, $\eta\xi^\pm_1$&  &  &  &  &  &  \\
    $\Nb\xi^\pm_2$, $\eta\xi^\pm_2$&  &  &  &  &  &\T\\
    $\Nb\xi_3$, $\eta\xi_3$        &  &  &  &  &\T&  \\
    $\Nb\xi_4$, $\eta\xi_4$        &  &  &  &\T&\T&\T\\
    \midrule
    $\xi^\pm_1\otimes\xi^\pm_1$    &\T&  &  &\T&  &  \\
    $\xi_2^\pm\otimes\xi_2^\pm$    &\T&\T&  &\T&\T&  \\
    $\xi_3\otimes\xi_3$            &\T&\T&  &\T&\T&\T\\
    $\xi_4\otimes\xi_4$            &\T&\T&  &\T&\T&\T\\
    \midrule
    $\xi^+_1\odot\xi^-_1$          &  &  &  &  &  &  \\
    $\xi^\pm_1\odot\xi^\pm_2$      &  &\T&  &  &\T&  \\
    $\xi^\pm_1\odot\xi^\mp_2$      &  &  &  &  &  &  \\
    $\xi^\pm_1\odot\xi_3$          &  &  &  &  &  &\T\\
    $\xi^\pm_1\odot\xi_4$          &  &  &\T&  &  &  \\
    $\xi^+_2\odot\xi_2^-$          &  &\T&  &  &\T&  \\
    $\xi^\pm_2\odot\xi_3$          &  &  &\T&  &  &\T\\
    $\xi^\pm_2\odot\xi_4$          &  &  &\T&  &  &\T\\
    $\xi_3\odot\xi_4$              &  &\T&  &  &\T&  \\
    \bottomrule
  \end{tabular}
  \caption{Ricci curvatures, $2n=6$}
  \label{tab:Ric6}
\end{table}

\begin{table}[tp]
  \centering
  \begin{tabular}{lcccccccc}
    \toprule
    &\multicolumn{2}{c}{\eqref{tresricmasric}}
    &\multicolumn{2}{c}{\eqref{richmenosrich}}
    &\eqref{ricastsh}
    &\eqref{ricsh}
    &\multicolumn{2}{c}{\cite{Falcitelli-FS:aH}}\\
    \cmidrule(rl){2-3}
    \cmidrule(rl){4-5}
    \cmidrule(rl){6-6}
    \cmidrule(rl){7-7}
    \cmidrule(rl){8-9}
    $2n=6$&$\Kc{1}$&$\Kc{2}$&$\Kc{-1}$&$\Kc{-2}$&$\Cc6
    $&$\Cc8$&$\Cc5$&$\Cc7$\\
    \midrule
    $d\Ieta$                       &\T&\T&  &  &  &  &  &  \\
    \midrule
    $\Nb\xi^\pm_1$, $\eta\xi^\pm_1$&  &  &  &  &\T&  &  &  \\
    $\Nb\xi^\pm_2$, $\eta\xi^\pm_2$&  &  &  &  &\T&\T&\T&\T\\
    $\Nb\xi_3$, $\eta\xi_3$        &  &\T&  &\T&  &  &  &\T\\
    $\Nb\xi_4$, $\eta\xi_4$        &\T&\T&\T&\T&\T&\T&  &  \\
    \midrule
    $\xi^\pm_1 \otimes \xi^\pm_1$  &\T&  &\T&  &  &  &  &  \\
    $\xi^\pm_2 \otimes \xi^\pm_2$  &\T&\T&\T&\T&  &  &  &  \\
    $\xi_3   \otimes \xi_3$        &\T&\T&  &  &  &\T&  &\T\\
    $\xi_4   \otimes \xi_4$        &\T&\T&  &  &  &\T&  &  \\
    \midrule
    $\xi^+_1 \odot \xi^-_1$        &  &  &  &  &  &  &  &  \\
    $\xi^\pm_1 \odot \xi^\pm_2$    &  &\T&  &\T&  &  &  &  \\
    $\xi^\pm_1 \odot \xi^\mp_2$    &  &  &  &  &  &  &  &  \\
    $\xi^\pm_1 \odot \xi_3$        &  &  &  &  &  &\T&\T&  \\
    $\xi^\pm_1 \odot \xi_4$        &  &  &  &  &\T&  &  &  \\
    $\xi^+_2 \odot \xi^-_2$        &  &\T&  &\T&  &  &  &  \\
    $\xi^\pm_2 \odot \xi_3$        &  &  &  &  &\T&\T&\T&\T\\
    $\xi^\pm_2 \odot \xi_4$        &  &  &  &  &\T&\T&\T&\T\\
    $\xi_3   \odot \xi_4$          &  &\T&  &  &  &  &  &\T\\
    \bottomrule
  \end{tabular}
  \caption{Curvature complementary to $\Cc3=\Cur(\su(n))$, $2n=6$}
  \label{tab:Cur6}
\end{table}

\subsection{Four dimensions}

The $\Un(2)$-decomposition of the space of curvature tensors $\FCur$ is
given by
\begin{equation*}
  \FCur = \Kah + \Kah^\perp = \Cc3 + \Kc1 + \Kc2 + \Kc{-1} +
  \Cc5 + \Cc6 + \Cc8.
\end{equation*}
When we consider the $\SU(2)$ action, only the modules $\Cc3$, $\Kc{1}$,
$\Kc2$ and $\Kc{-1}$ remain irreducible.  To describe the decompositions of
$\Cc5$ and $\Cc6$ into $\SU(2)$-irreducible modules, we will make use of
tensors defined by
\begin{equation*}
  \chi(a,b) = 6\, a \odot b - a \wedge b,
\end{equation*}
for all $a,b \in \Lambda^2 T^* M$, where $\odot$ denotes the symmetric
product given by $2\,a\odot b = a \otimes b + b \otimes a$.  The relevant
decompositions are now given by
\begin{enumerate}
\item $\Cc5 = \Cc5^{++} + \Cc5^{--} + \Cc5^{+-}$, where $\Cc5^{++} =
  \mathbb R \chi(\psi_+, \psi_+)$, $\Cc5^{--} = \mathbb R\chi(\psi_- ,
  \psi_-)$ and $\Cc5^{+-} = \mathbb R \chi(\psi_+ , \psi_-)$,
\item $ \Cc6 = \Cc6^{+} + \Cc6^{-} $, where $ \Cc6^{+} = \mathbb R
  \chi(\psi_+ , \omega)$ and $\Cc6^{-} = \mathbb R \chi(\psi_- , \omega )$.
\end{enumerate}
For the intrinsic torsion, the $\Un(2)$-decomposition of $\xi$ is given by
\begin{equation*}
\xi = \xi_2 + \xi_4 \in \mathcal W = \Wc2 + \Wc4.
\end{equation*}
Under $\SU(2)$, we have $\Wc2 \cong \Wc4 \cong T^* M$, which we will see
gives rise to different choices of decompositions of~$\xi$.

For an $\SU(2)$-structure, we have $\nabla \omega \in \mathcal W = T^* M
\otimes \psi_+ + T^* M \otimes \psi_-$.  Consequently, $\nabla \omega =
\xi_+ \otimes \psi_+ + \xi_- \otimes \psi_-$, where $\xi_+$ and $\xi_-$ are
one-forms.  Moreover,
\begin{equation*}
2 \langle Y , \xi_X Z \rangle = - \xi_+(X) \psi_-(Y,Z) + \xi_-(X)
\psi_+(Y,Z),
\end{equation*}
so $\xi = \xi_+ + \xi_-$, where
\begin{equation*}
  2\inp Y{(\xi_+)_XZ} = - \xi_+(X) \psi_-(Y,Z),\qquad
  2\inp Y{(\xi_-)_XZ} = \xi_-(X) \psi_+(Y,Z).
\end{equation*}

The two decompositions of $\xi$ are related as follows:
\begin{enumerate}
\item $\xi \in \Wc2$ if and only if $\xi_+ = I \xi_-$.
\item $\xi \in \Wc4$ if and only if $\xi_+ = - I \xi_-$.
\end{enumerate}

The following theorem gives information about the contributions of the
components of the intrinsic torsion to the tensors $\Ric^*$ and $\Ric$.  We
first note that in dimension four, $\SRic^*$~decomposes under $\SU(2)$ as
\begin{equation*}
  \SRic^* = \mathbb R\inp\cdot\cdot + [\lambda_0^{1,1}] +
  \mathbb  R \psi_+ + \mathbb R \psi_-.
\end{equation*}

\begin{theorem}
  Let $M$ be a special almost Hermitian $4$-manifold with minimal
  $\SU(2)$-connection $\Nb = \nabla + \eta + \xi = \Nt + \eta $.  The
  curvature contributions corresponding to Theorems~\ref{thm:Ric}
  and~\ref{thm:Curv} via the decompositions $\xi=\xi_2+\xi_4$ and
  $\xi=\xi_++\xi_-$ are given in Tables \ref{tab:Ric4} and~\ref{tab:Cur4}.
\end{theorem}

\begin{table}[tp]
  \centering
  \begin{tabular}{lccccccc}
    \toprule
    &\multicolumn{4}{c}{$\Ric^*$
    \eqref{ricasteta}}&\multicolumn{3}{c}{$\Ric$ \eqref{riceta}}\\
    \cmidrule(lr){2-5}
    \cmidrule(lr){6-8}
    $2n=4$&$\mathbb R$&$[\lambda^{1,1}_0]$&$\mathbb
    R\psi_+$&$\mathbb R\psi_-$&$\mathbb
    R$&$[\lambda^{1,1}_0]$&$\real{\sigma^{2,0}}$\\  
    \midrule
    $d\Ieta$                   &\T&\T&\T&\T&\T&\T&  \\
    \midrule                         
    $\Nb\xi_2$, $\eta\xi_2$    &  &  &  &  &  &  &\T\\
    $\Nb\xi_4$, $\eta\xi_4$    &  &  &  &  &\T&  &\T\\
    $\xi_2\otimes\xi_2$        &\T&\T&  &  &\T&\T&  \\
    $\xi_4\otimes\xi_4$        &\T&\T&  &  &\T&\T&\T\\
    $\xi_2\odot\xi_4$          &  &  &  &  &  &\T&\T\\
    \midrule                         
    $\Nb\xi_+$, $\eta\xi_+$    &  &  &  &  &\T&  &\T\\
    $\Nb\xi_-$, $\eta\xi_-$    &  &  &  &  &\T&  &\T\\
    $\xi_+\otimes\xi_+$        &  &  &  &  &\T&  &\T\\
    $\xi_-\otimes\xi_-$        &  &  &  &  &\T&  &\T\\
    $\xi_+\odot\xi_-$          &\T&\T&  &  &\T&\T&\T\\
    \bottomrule
  \end{tabular}
  \caption{Ricci curvatures, $2n=4$}
  \label{tab:Ric4}
\end{table}

\begin{table}[tp]
  \centering
  \begin{tabular}{lccccccccc}
    \toprule
    $2n=4$&$\Kc1$&$\Kc2$&$\Kc{-1}$& $\Cc6^+$&$\Cc6^-$&$\Cc8$&$\Cc5^{++}$
    &$\Cc5^{--}$&$\Cc5^{+-}$\\
    \midrule
    $d\Ieta$                &\T&\T&  &  &  &  &  &  &  \\
    \midrule
    $\Nb\xi_2$, $\eta\xi_2$ &  &  &  &\T&\T&\T&\T&\T&\T\\
    $\Nb\xi_4$, $\eta\xi_4$ &\T&  &\T&\T&\T&\T&  &  &  \\
    $\xi_2 \otimes \xi_2 $  &\T&\T&\T&  &  &  &  &  &  \\
    $\xi_4 \otimes \xi_4 $  &\T&\T&  &  &  &\T&  &  &  \\
    $\xi_2 \odot \xi_4 $    &  &  &  &\T&\T&\T&\T&\T&\T\\
    \midrule
    $\Nb\xi_+$, $\eta\xi_+$ &\T&  &\T&  &\T&\T&\T&\T&\T\\
    $\Nb\xi_-$, $\eta\xi_- $&\T&  &\T&\T&  &\T&\T&\T&\T\\
    $\xi_+ \otimes  \xi_+ $ &\T&  &\T&  &  &\T&  &\T&\T\\
    $\xi_- \otimes  \xi_- $ &\T&  &\T&  &  &\T&\T&  &\T\\
    $\xi_+ \odot  \xi_- $   &\T&\T&\T&\T&\T&\T&\T&\T&\T\\
    \bottomrule
  \end{tabular}
  \caption{Curvature complementary to $\Cc3=\Cur(\su(n))$, $2n=4$}
  \label{tab:Cur4}
\end{table}

\begin{proof}
  The absence of~$\Kc{-2}$ in the decomposition of~$\FCur$ comes from the
  fact that
  \begin{equation} \label{escalar}
    \left(\Ric^*_H - \Ric_H   \right) (X,Y) = \beta\, \inp XY,
  \end{equation}
  where $\beta = \inp{(\Nt_{e_i}\xi)_{e_j}e_j}{e_i} +
  \inp{\xi_{\xi_{e_i}e_j}e_j}{e_i}$.  Therefore, by \eqref{ricastsh}, we
  have
  \begin{equation}
    \label{k1k2}
    \begin{split}
      \left(3 \Ric^*_H + \Ric_H \right) (X,Y)
      & =  - \beta\inp XY - 4 d\Ieta(X,IY) + 4 d\Ieta(IX,Y)\\
      &\qquad - 2 \inp{\xi_Xe_i}{\xi_{IY}Ie_i} - 2
      \inp{\xi_Ye_i}{\xi_{IX}Ie_i}. 
    \end{split}
  \end{equation}
  Using equations \eqref{escalar} and \eqref{k1k2} the tables follow.
\end{proof}

\begin{remark}
  Let us list some direct consequences of results and tables presented here
  and in~\S\ref{sec:high}:
  \begin{enumerate}[(a)]
  \item if $\xi \in \Wc3$, then the components of $R$ in $\Kc{-1}$, $\Cc5$
    and $\Cc6$ vanish;
  \item if $\xi\in \Wc3 + \Wc4$ and $d\Ieta$ is Hermitian,
    then the components of $R$ in $\Cc5$ and $\Cc6$ vanish;
  \item if $\xi \in \Wc1 + \Wc2$ and $d\Ieta$ is Hermitian,
    then the component of $R$ in $\Cc6$ vanishes; and
  \item if $n=2$ and $d\Ieta$ is Hermitian, then the component of~$R$
    in~$\Cc6$ vanishes.
  \end{enumerate}
  There are more consequences of this sort, but they have been already
  pointed out in~\cite{Falcitelli-FS:aH}.
\end{remark}

\begin{remark}
  For special almost Hermitian $2$-manifolds, we have the following
  identity, deduced in \cite{Cabrera:special},
  \begin{equation*}
    K(\psi_+ , \psi_-) = d\Ieta (\psi_+ , \psi_-) = d \eta_+ (\psi_+) +
    d \eta_- (\psi_-) - \eta_+^2 - \eta_-^2,
  \end{equation*}
  where $K$ denotes the sectional curvature and $\Ieta = \eta_+ \psi_- -
  \eta_- \psi_+$.
\end{remark}

\section{Identities from the exterior algebra}
\label{sec:exterior}

As remarked in~\S\ref{sec:high} one may see different contributions to the
module $\Cc6\cong\real{\lambda^{2,0}}$ by using different computations of
the curvature.  This is because of non-trivial identities relating the
components of $\Nt\xi_i$ and $\xi_j\odot\xi_k$.  Such an identity for the
$\real{\lambda^{2,0}}$-components may be obtained by comparing equations
\eqref{ricastsh} and~\eqref{otraricsh}.  However, we claim that this
information may also be obtained from the exterior algebra of a
$\Un(n)$-manifold.

Consider the K\"ahler two-form~$\omega$.  Being a differential form it
satisfies $d^2\omega=0$.  However, since the Levi-Civita
connection~$\nabla$ is torsion-free, we may compute $d^2\omega$ using
$\nabla$.  Writing $\nabla=\Nt-\xi$ and using $\Nt\omega=0$, we have first
that
\begin{equation*}
  \tfrac12d\omega(Y,Z,W) = \inp{\xi_YZ}{IW} + \inp{\xi_WY}{IZ} +
  \inp{\xi_ZW}{IY}. 
\end{equation*}
Now $d^2\omega = \alt(\Nt d\omega)- \alt(\xi d\omega)$, where $\alt\colon
T^*M\otimes \Lambda^3T^*M\to\Lambda^4T^*M$~is the alternation map.  One
computes that these two terms are the expressions obtained respectively by
summing $\varepsilon \inp{(\Nt_X\xi)_YZ}{IW}$ and $\varepsilon
\inp{\xi_{\xi_XY}Z}{IW}$ over all permutations of $(X,Y,Z,W)$, where
$\varepsilon$~is the sign of the permutation.

We have that 
\begin{equation*}
  \Lambda^4T^*M =
  \real{\lambda^{4,0}}+\real{\lambda^{3,1}}+\real{\lambda^{2,0}}\omega +
  [\lambda^{2,2}_0]+[\lambda^{1,1}_0]\omega+\mathbb R\omega^2,
\end{equation*}
so in order to compute the $\real{\lambda^{2,0}}$-component of $d^2\omega$
we contract with $\omega$ on the first two arguments and then take the
projection to $\real{\lambda^{2,0}}$, which is the $(-1)$-eigenspace of $I$
acting on $2$-forms.  Using the symmetries of the components of~$\xi$ one
obtains that the $\real{\lambda^{2,0}}$-component of $d^2\omega$ is
\begin{equation}
  \label{d2omega:part20}
  \begin{split}
    0
    &=3\inp{(\Nt_{e_i}\xi_1)_{e_i}X}Y - \inp{(\Nt_{e_i}\xi_3)_{e_i}X}Y
    + (n-2) \inp{(\Nt_{e_i}\xi_4)_{e_i}X}Y \\
    &\quad
    + \inp{{\xi_3}_Xe_i}{{\xi_1}_{e_i}Y}
    - \inp{{\xi_3}_Ye_i}{{\xi_1}_{e_i}X}
    + \inp{{\xi_3}_Xe_i}{{\xi_2}_{e_i}Y}
    - \inp{{\xi_3}_Ye_i}{{\xi_2}_{e_i}X}\\
    &\quad
    - \frac{n-5}{n-1}\inp{{\xi_1}_{{\xi_4}_{e_i}e_i}X}Y
    - \frac{n-2}{n-1}\inp{{\xi_2}_{{\xi_4}_{e_i}e_i}X}Y
    + \inp{{\xi_3}_{{\xi_4}_{e_i}e_i}X}Y.
  \end{split}
\end{equation}
We conclude that in general dimensions there is a non-trivial linear
relation between the $\real{\lambda^{2,0}}$-components of $\Nt\xi_1$,
$\Nt\xi_3$, $\Nt\xi_4$, $\xi_1\odot\xi_3$, $\xi_1\odot\xi_4$,
$\xi_2\odot\xi_3$, $\xi_2\odot\xi_4$ and $\xi_3\odot\xi_4$.  By
`non-trivial' we mean that no coefficient is zero, so this relation may be
used to write any of the terms as a linear combination of the others.
Interestingly, this relation does not involve $\xi_1\odot\xi_4$, when
$2n=10$.

This is sufficient to explain the difference between the ticks in the
$\Cc6$ column in~\cite{Falcitelli-FS:aH} and those we obtained from
equation~\eqref{ricastsh}.  An extra coincidence in the coefficients
explains the differences between our results from~\eqref{ricastsh}
and~\eqref{otraricsh}.

One may try to apply the above approach to the other modules that
$\Lambda^4T^*M$ has in common with the space of curvature tensors, namely
$[\lambda^{2,2}_0]$, $[\lambda^{1,1}_0]\omega$ and $\mathbb R\omega^2$.
However, this is not so rewarding because of the higher multiplicities that
these modules have in the relevant decompositions.  Indeed,
$\Cc6\cong\real{\lambda^{2,0}}$ is distinguished by occurring only with
multiplicity one or zero in the modules for $\Nt\xi_i$
and~$\xi_i\otimes\xi_j$.

In~\cite{Falcitelli-FS:aH}, it is pointed out that if $\xi \in \Wc4$, then
the components of $R$ in $\Cc4$, $\Cc5$, $\Cc6$ and $\Cc7$ vanish.  Let us
indicate how equation~\eqref{d2omega:part20} gives an alternative proof of
this result, for $n>2$. In fact, using Tables \ref{tab:Cur8}
and~\ref{tab:Cur6}, the vanishing of the components in $\Cc4$, $\Cc5$
and~$\Cc7$ is immediate.  On the other hand, equations \eqref{ricastsh}
and~\eqref{d2omega:part20} give the vanishing of the component in~$\Cc6$.

Finally, a comparison of Tables \ref{tab:Ric8} and~\ref{tab:Cur8} reveals
another relation on special almost Hermitian manifolds: the
$\real{\lambda^{2,0}}$-part of $d\Ieta$ carries all the information from
the corresponding components of $\Nt\xi_i$ modulo the
$\real{\lambda^{2,0}}$-parts of $\xi_1\odot\xi_3$, $\xi_1\odot\xi_4$,
$\xi_2\odot\xi_3$ and $\xi_2\odot\xi_4$.  This relation is obtainable by
considering the $(n,2)$-part of the equation $d^2\Psi=0$, where $\Psi$~is
the complex volume, cf.~\cite{Cabrera:special}.

\begin{small}

\begin{thebibliography}{10}

\bibitem{Banos:MA6}
B.~Banos, \emph{Nondegenerate {M}onge-{A}mp{\`e}re structures in
  dimension~{$6$}}, Lett. Math. Phys. \textbf{62} (2002), no.~1, 1--15.

\bibitem{Besse:Einstein}
A.~L. Besse, \emph{{E}instein manifolds}, Ergebnisse der Mathematik und ihrer
  Grenzgebiete, 3. Folge, vol.~10, Springer, Berlin, Heidelberg and New York,
  1987.

\bibitem{Bryant:splag}
R.~L. Bryant, \emph{Some examples of special {L}agrangian tori}, Adv. Theor.
  Math. Phys. \textbf{3} (1999), no.~1, 83--90.

\bibitem{Chiossi-S:SU3-G2}
S.~G. Chiossi and S.~Salamon, \emph{The intrinsic torsion of {$\rm SU(3)$} and
  {$G\sb 2$} structures}, Differential geometry, {V}alencia, 2001, World Sci.
  Publishing, River Edge, NJ, 2002, pp.~115--133.

\bibitem{Cleyton-S:intrinsic}
R.~Cleyton and A.~F. Swann, \emph{{E}instein metrics via intrinsic or parallel
  torsion}, Math. Z. \textbf{247} (2004), no.~3, 513--528.

\bibitem{Falcitelli-FS:aH}
M.~Falcitelli, A.~Farinola, and S.~M. Salamon, \emph{Almost-{H}ermitian
  geometry}, Differential Geom. Appl. \textbf{4} (1994), 259--282.

\bibitem{Grantcharov-GP:CY-toric}
D.~Grantcharov, G.~Grantcharov, and Y.~S. Poon, \emph{{C}alabi-{Y}au
  connections with torsion on toric bundles}, June 2003, eprint
  \url{arXiv:math.DG/0306207}.

\bibitem{Gray:curvature}
A.~Gray, \emph{Curvature identities for {H}ermitian and almost {H}ermitian
  manifolds}, T\^ohoku Math. J. (2) \textbf{28} (1976), no.~4, 601--612.

\bibitem{Gray:nearly-Kaehler}
\bysame, \emph{The structure of nearly {K\"a}hler manifolds}, Math. Ann.
  \textbf{223} (1976), 233--248.

\bibitem{Gray-H:16}
A.~Gray and L.~M. Hervella, \emph{The sixteen classes of almost {H}ermitian
  manifolds and their linear invariants}, Ann. Mat. Pura Appl. (4) \textbf{123}
  (1980), 35--58.

\bibitem{Gutowski-IP:calibrations}
J.~Gutowski, S.~Ivanov, and G.~Papadopoulos, \emph{Deformations of generalized
  calibrations and compact non-{K}\"ahler manifolds with vanishing first
  {C}hern class}, Asian J. Math. \textbf{7} (2003), no.~1, 39--79.

\bibitem{Hitchin:special-L}
N.~J. Hitchin, \emph{The moduli space of special {L}agrangian submanifolds},
  Ann. Scuola Norm. Sup. Pisa Cl. Sci. (4) \textbf{25} (1997), no.~3-4,
  503--515 (1998), Dedicated to Ennio De Giorgi.

\bibitem{Joyce:holonomy}
D.~Joyce, \emph{Compact manifolds with special holonomy}, Oxford Mathematical
  Monographs, Oxford University Press, 2000.

\bibitem{Cabrera:special}
F.~Mart{\'\i}n~Cabrera, \emph{Special almost {H}ermitian geometry}, eprint
  \url{arXiv:math.DG/0409167}, September 2004.

\bibitem{Papadopoulos:brane}
G.~Papadopoulos, \emph{Brane solitons and hypercomplex structures},
  Quaternionic structures in mathematics and physics (Rome, 1999), Univ. Studi
  Roma ``La Sapienza'', Rome, 1999, pp.~299--312 (electronic).

\bibitem{Salamon:holonomy}
S.~M. Salamon, \emph{Riemannian geometry and holonomy groups}, Pitman Research
  Notes in Mathematics, vol. 201, Longman, Harlow, 1989.

\bibitem{Tricerri-Vanhecke:aH}
F.~Tricerri and L.~Vanhecke, \emph{Curvature tensors on almost {H}ermitian
  manifolds}, Trans. Amer. Math. Soc. \textbf{267} (1981), 365--398.

\end{thebibliography}
\providecommand{\bysame}{\leavevmode\hbox to3em{\hrulefill}\thinspace}
\providecommand{\MR}{\relax\ifhmode\unskip\space\fi MR }
\providecommand{\MRhref}[2]{%
  \href{http://www.ams.org/mathscinet-getitem?mr=#1}{#2}
}
\providecommand{\href}[2]{#2}

\end{small}

\begin{smallpars}
  Mart\'\i n Cabrera: \textit{Department of Fundamental Mathematics,
  University of La Laguna, 38200 La Laguna, Tenerife, Spain}.  E-mail:
  \url{fmartin@ull.es}
  
  Swann: \textit{Department of Mathematics and Computer Science, University
  of Southern Denmark, Campusvej 55, DK-5230 Odense M, Denmark}.  E-mail:
  \url{swann@imada.sdu.dk}
\end{smallpars}

\end{document}